\documentclass[english]{article}

\usepackage{amsmath,amsthm,amssymb,babel}

\def\qed{\hfill $\Box$\par\vskip3mm}
\def\di{\displaystyle}
\def\ri{\rightarrow}
\def\ep{\varepsilon}
\def\RR{{\mathbb R}}
\def\ee{{\mathcal E}}

\def\plapl{{\rm div}\,(|\nabla u|^{p-2}\nabla u)}
\def\d1p{{\mathcal D}_0^{1,p}(\Omega)}
\def\d1pr{{\mathcal D}^{1,p}(\RR^N)}
\def\intom{\int_\Omega}
\def\intor{\int_{\RR^N}}
\newtheorem{theorem}{Theorem}[section]
\newtheorem{lemma}{Lemma}

\title{\sc Multiple perturbations of a singular eigenvalue problem}

\author{Matija Cencelj,  Du\v san Repov\v s and \v Ziga Virk}  
\date{}

\begin{document}
\maketitle

\noindent{\small {\bf Abstract}.  We study the perturbation by a critical term and a $(p-1)$-superlinear subcritical nonlinearity of a quasilinear elliptic equation containing a singular potential. By means of variational arguments and a version of the concentration-compactness principle in the singular case, we prove the existence of solutions for positive values of the parameter under the principal eigenvalue of the associated singular eigenvalue problem.\\
{\bf Key words}: critical Sobolev exponent; indefinite singular potential; principal eigenvalue; concentration-compactness; singularity.\\
{\bf 2010 Mathematics Subject Classification}: 35B33; 35B38; 35P30; 47J10; 58E05.}

\section{Introduction}
Let $\Omega\subseteq\RR$ be an {\it arbitrary} open set, $1<p<N$, and let ${\mathcal D}_0^{1,p}(\Omega)$ denote the completion of ${\mathcal D}(\Omega)$ with respect to the norm
$\|u\|:=(\intom |\nabla u|^pdx)^{1/p}$. Let $V\in L^1_{loc}(\Omega)$ be a function which may have strong singularities and an indefinite sign.

 Smets was interested in \cite{smets} in finding nontrivial weak solutions for the following nonlinear eigenvalue problem:
\begin{equation}\label{psmets}
\left\{\begin{array}{ll}
&\di -\plapl =\lambda\,V(x)|u|^{p-2}u\quad\mbox{in}\ \Omega\\
&\di u\in{\mathcal D}_0^{1,p}(\Omega)\,.\end{array}\right.
\end{equation}

Problems of this type are in relationship with the study of the standing waves in anisotropic
Schr\"odinger or Klein-Gordon equations, cf. Reed and Simon \cite{reed}, Strauss \cite{str}, and Wang \cite{wang}. Equation \eqref{psmets} is also considered
 a model for several physical phenomena related to the equilibrium of anisotropic
media that possibly are somewhere
{\it perfect insulators}
or {\it perfect conductors}, see Dautray and Lions \cite[p. 79]{dautray}. We point out that degenerate or singular problems have been intensively studied
starting with the pioneering paper by Murthy and Stampacchia \cite{murthy}.

Problem \eqref{psmets} is in relationship with several papers dealing with nonlinear anisotropic eigenvalue problems, see Brown and Tertikas \cite{brown}, Rozenblioum and Solomyak \cite{roze}.  Szulkin and Willem generalize in \cite{szulkin} several earlier results concerning the existence of an infinite
sequence of eigenvalues. The main hypothesis on the potential $V$ in \cite{szulkin} is the following:
\begin{equation}\label{hypsw}\left\{\begin{array}{ll}
&\di V\in L^1_{loc}(\Omega),\ V^+=V_1+V_2\not=0,\ V_1\in L^{N/p}(\Omega),\\
&\di \mbox{for every}\ y\in\overline\Omega,\ \lim_{x\ri y,x\in\Omega}|x-y|^pV_2(x)=0\ \mbox{and}\\
&\di \lim_{x\ri\infty ,x\in\Omega}|x|^pV_2(x)=0.
\end{array}\right.
\end{equation}

Under  assumption \eqref{hypsw}, the mapping ${\mathcal D}_0^{1,p}(\Omega)\ni u\longmapsto\intom V^+|u|^pdx$ is weakly continuous, so the problem is not affected by a lack of compactness.
In \cite{smets}  the case of indefinite potential functions $V$ is studied
for which no {\it a priori} compactness is assumed. The corresponding hypotheses extend condition \eqref{hypsw}, nonetheless they are not directly linked to punctual growths of $V$.
Due to the presence of a singular potential, the classical methods cannot be applied directly, so the existence
can become
a delicate matter.

Consider the minimization problem
\begin{equation}\label{mV}
S_V:=\inf\left\{\intom|\nabla u|^pdx;\ u\in{\mathcal D}_0^{1,p}(\Omega),\ \intom V(x)|u|^pdx=1\right\}.\end{equation}
As established in \cite{smets} with standard constrained minimization arguments, minimizers of problem \eqref{mV} correspond to weak solutions of \eqref{psmets}, with $\lambda$ appearing as a Lagrange multiplier (that is, $\lambda =S_V$).
Such a parameter $\lambda$ is called the {\it principal eigenvalue} for problem \eqref{psmets}.

In order to have $S_V\not= 0$ and well defined, we assume that $V=V^+-V_-$, $V^+\not=0$, and that there exists $c>0$ such that for all $u\in{\mathcal D}_0^{1,p}(\Omega)$,
\begin{equation}\label{h1}
c\intom V^+|u|^pdx\leq\intom|\nabla u|^pdx.\end{equation}
By Hardy's inequality it follows that potentials with point singularities and decay at infinity both at most as $ O(|x|^{-p})$ satisfy hypothesis \eqref{h1}.

Since $\Omega$ is not necessarily bounded and $V$ can have singularities, it is not clear that the infimum in problem \eqref{mV} is achieved without imposing additional conditions that allow the analysis of minimizing sequences. For all $x\in\overline\Omega$ and $r>0$, we denote by $B_r(x)$  the open  ball centered at $x$ and of radius $r$ and by $B_r$ the closed ball centered at the origin (we can assume without any loss of generality that $0\in\Omega$). We introduce the following quantities:
$$S_{r,V}:=\inf\left\{\intom|\nabla u|^pdx;\ u\in {\mathcal D}(\Omega\setminus B_r),\ \intom V^+(x)|u|^pdx=1\right\};$$
$$S_{\infty ,V}:=\sup_{r>0}S_{r,V}=\lim_{r\ri\infty}S_{r,V};$$
$$S_{r,V}^x:=\inf\left\{\intom|\nabla u|^pdx;\ u\in {\mathcal D}(\Omega\cap B_r(x)),\ \intom V^+(x)|u|^pdx=1\right\};$$
$$S^x_V:=\sup_{r>0}S^x_{r,V}=\lim_{r\ri 0}S^x_{r,V};$$
$$S_{*,V}:=\inf_{x\in\overline\Omega}S^x_V;$$
$$\Sigma_V:=\{x\in\overline\Omega;\ S^x_V<\infty\}.$$

Applying Hardy's inequality
$$\int_{\RR^N}\frac{|u|^p}{|x|^p}\,dx\leq\left(\frac{N}{N-p}\right)^p\int_{\RR^N}|\nabla u|^pdx,$$
we observe that under assumption \eqref{hypsw} introduced in \cite{szulkin}, we have $S_{\infty ,V}=S_{*,V}=+\infty$. As argued in \cite[p. 475]{smets}, the condition $S_{\infty ,V}=S_{*,V}=+\infty$ is equivalent to the weak continuity of the mapping $u\longmapsto\intom V^+(x)|u|^pdx$.

We make the following hypothesis:
\begin{equation}\label{h2}
\mbox{the closure of $\Sigma_V$ is at most countable.}\end{equation}

In particular, condition \eqref{h2} excludes the presence of strong {\it spikes} on a dense subset of $\Omega$.

 For $V\in L^1_{loc}(\Omega)$ satisfying assumptions \eqref{h1} and \eqref{h2}, Smets proved in \cite{smets} that the singular eigenvalue problem \eqref{psmets} admits a principal eigenvalue, provided that $S_V<S_{\infty,V}$ and $S_V<S_{*,V}$. This result extends and simplifies the work of Tertikas \cite{tertikas}, which deals with the positive linear case for $\Omega=\RR^N$.
We point out (see \cite[p. 472]{smets}) that the condition $p<N$ is necessary only if $\Omega$ is unbounded, otherwise one can work in the standard Sobolev space $W_0^{1,p}(\Omega)$.

We are interested in studying what happens if problem \eqref{psmets} is affected by  certain perturbations. This is needed in several applications
and the idea of using perturbation methods in the
treatment of nonlinear boundary value problems was introduced by Struwe \cite{struwe}.
Existence results for nonautonomous perturbations of critical singular elliptic boundary value problems were established by R\u adulescu and Smets \cite{rasm}; in their case, the singular weight allows for
unbounded domains as cones and gives rise to a different noncompactness picture, as was first
remarked by Caldiroli and Musina \cite{camu}.

Let ${\mathcal M}(\RR^N)$ denote the Banach space of finite Radon measures over $\RR^N$ endowed with the norm
$$\|\mu\|:=\sup_{\phi\in C_0(\RR^N),|\phi|_\infty\leq 1}|\mu(\phi)|.$$
By definition, a sequence $(\mu_n)\subset {\mathcal M}(\RR^N)$ weakly converges  to $\mu\in {\mathcal M}(\RR^N)$ if $\mu_n(\phi)\ri\mu(\phi)$ for all $\phi\in C_0(\RR^N)$. The Banach-Alaoglu theorem implies that every bounded sequence $(\mu_n)\subset {\mathcal M}(\RR^N)$ contains a weakly convergent subsequence. We denote by ${\mathcal M}^+(\RR^N)$ the cone of positive Radon measures over $\RR^N$ and by $\delta_x$ the Dirac mass at the point $x$.

\section{Effects of a double perturbation}
In the present paper, we are concerned with a perturbation of problem \eqref{psmets} and we are interested in the combined effects of a $(p-1)$-superlinear subcritical nonlinearity and a critical Sobolev term. To fix the ideas, we consider $\Omega=\RR^N$ but the arguments can be adapted to any open set in $\RR^N$. More precisely, we study the nonlinear problem
\begin{equation}\label{problem}
\left\{\begin{array}{ll}
&\di -\plapl =\lambda\,V(x)|u|^{p-2}u+a(x)|u|^{r-2}u+b(x)|u|^{p^*-2}u\quad\mbox{in}\ \RR^N\\
&\di u\in\d1pr\,,\end{array}\right.
\end{equation}
where $p^*=Np/(N-p)$ stands for the critical Sobolev exponent.

This problem can be viewed as a prototype of pattern formation in biology
and is related to the steady-state problem for a chemotactic aggregation model
introduced by Keller and Segel \cite{keller}. Problem \eqref{problem} also plays a crucial role in
the analysis of activator-inhibitor systems modeling biological pattern formation, cf. Gierer and Meihardt \cite{gierer}.

Problem \eqref{problem} is related to the Brezis-Nirenberg problem
\begin{equation}\label{bn}
-\Delta u=\lambda u+u^{(N+2)/(N-2)}\quad\mbox{in}\ \Omega\subset\RR^N,\end{equation}
where $\Omega$ is an open bounded set with smooth boundary.
Brezis and Nirenberg \cite{brenir} showed that, contrary to intuition, the critical problem with {\it small} linear perturbation can provide solutions.
More precisely, Brezis and Nirenberg proved that
problem \eqref{bn} admits a positive solution vanishing on $\partial\Omega$  if and only if $0<\lambda<\lambda_1$ (if $N\geq 4$), where
$\lambda_1$ is the first eigenvalue of the Laplace operator in $H^1_0(\Omega)$. In \cite{brenir}, other results are also established (for instance, if $N=3$ or when $\lambda$ is replaced by $g(x,u)$ satisfying an appropriate growth condition) and pioneering techniques in nonlinear analysis  are introduced.

Our assumptions are the following:
\begin{equation}\label{hypr}
p<r<p^*;
\end{equation}
\begin{equation}\label{hypa}
a\in L^s(\RR^N)\ \mbox{with}\ s=\frac{Np}{Np-r(N-p)},\ a(x)\geq 0\ \mbox{a.e.}\ x\in\RR^N,\ a\not= 0;
\end{equation}
\begin{equation}\label{hypb}
b\in L^\infty(\RR^N),\ b(0)=\|b\|_{L^\infty(\RR^N)},\ b(x)=b(0)+o(|x|^\eta)\ \mbox{as}\ x\ri 0,\end{equation}
where
$$\eta=\frac{N(s-1)}{(p-1)s}\quad\mbox{if}\ N<\frac{pr}{r+1-p};$$
$$\eta=\frac{N}{s}\quad\mbox{if}\ N\geq\frac{pr}{r+1-p}.$$

The asymptotic decay of the potential $b$ described in condition \eqref{hypb} compensates for the critical behaviour of the corresponding nonlinearity and it provides a sufficient condition for the existence of the ``valley" in the mountain pass theorem.

The solutions of problem \eqref{problem} correspond to nontrivial critical points of the  energy functional $\ee:\d1pr\ri\RR$ defined by
$$\ee (u)=\frac 1p\intor|\nabla u|^pdx-\frac\lambda p\intor V(x)|u|^pdx-\frac 1r\intor a(x)|u|^rdx-\frac{1}{p^*}\intor b(x)|u|^{p^*}dx.$$

Let $\lambda_1$ denote the principal eigenvalue of problem \eqref{psmets}, namely $\lambda_1=S_V$ in the minimization problem \eqref{mV}.
As remarked in \cite[p. 464]{smets}, hypothesis \eqref{h1} implies that $\lambda_1>0$.
Our main result asserts that the perturbed problem \eqref{problem} admits nontrivial solutions for all positive parameters $\lambda$ less than the principal eigenvalue of problem \eqref{psmets}.

\begin{theorem}\label{t1}
Let $V\in L^1_{loc}(\RR^N)$ satisfy
 $S_V<S_{\infty,V}$, $S_V<S_{*,V}$, and hypotheses \eqref{h1}, \eqref{h2}. Assume that conditions \eqref{hypr}, \eqref{hypa}, and \eqref{hypb} are fulfilled. Then problem \eqref{problem}
 admits at least one nontrivial solution
 for all positive  parameters with $\lambda<\lambda_1$.
\end{theorem}

For $c\in\RR$, we recall that $\ee$ satisfies the localized Palais-Smale (PS)$_c$-condition if every sequence $(u_n)\subset\d1pr$ with $\ee (u_n)\ri c$ and $\ee'(u_n)\ri 0$ in $(\d1pr)'$, has a convergent subsequence in $\d1pr$.

The main idea of the proof of Theorem \ref{t1} is to apply the mountain pass theorem. Note that $p^*$ is the limiting Sobolev exponent for the embedding $\d1pr\subset L^{p^*}(\RR^N)$. Since this embedding is not compact, the functional $\ee$ does not satisfy the Palais-Smale condition. By using the $V$-dependent concentration-compactness principle of Smets \cite[Lemma 2.1]{smets}, we show that $\ee$ satisfies the localized (PS)$_c$-condition for certain values of $c$. In the final part of the proof, we argue that the geometric hypotheses of the mountain pass theorem are also fulfilled.

\section{The localized Palais-Smale condition}
In this section we assume that the hypotheses of Theorem \ref{t1} are satisfied and we are interested to find a range of values for $c>0$ such that $\ee$ satisfies the Palais-Smale (PS)$_c$-condition. An important role in this choice of $c$ is played by the Sobolev constant
\begin{equation}\label{6}
S:=\inf\left\{\intor |\nabla u|^pdx;\ u\in W^{1,p}(\RR^N),\ \intor |u|^{p^*}dx=1\right\}.\end{equation}
This corresponds to the best constant for the Sobolev embedding $W^{1,p}(\RR^N)\subset L^{p^*}(\RR^N)$. We recall (see Brezis and Nirenberg \cite[p. 443]{brenir}) some basic properties of this constant:

(i) $S$ can be defined for any open set $\Omega$, is independent of $\Omega$, and depends only on $N$.

(ii) The infimum in \eqref{6} is never achieved in the case of bounded open sets.

(iii) For the whole Euclidean space, the infimum in \eqref{6} is achieved by the function
\begin{equation}\label{ueps}u_\varepsilon (x)=C_\varepsilon\left(\varepsilon^{p/(p-1)}+|x|^{p/(p-1)}\right)^{-\frac{N-p}{p}},\end{equation}
for all $\varepsilon>0$, where $C_\varepsilon$ is a positive constant depending on $\varepsilon$.

Let $(u_n)\subset\d1pr$ be such that $\ee (u_n)\ri c$ and $\ee'(u_n)\ri 0$ in $(\d1pr)'$.
We find an interval $(0,c_0)$ such that $(u_n)$ contains a convergent subsequence, provided that $c\in (0,c_0)$.
For this purpose we use some ideas found in the  paper by Guedda and V\'eron \cite{guedda}. We have
\begin{equation}\label{1}\begin{array}{ll}
\di \frac 1p\intor|\nabla u_n|^pdx&\di -\frac\lambda p\intor V(x)|u_n|^pdx-\frac 1r\intor a(x)|u_n|^rdx\\
&\di -\frac{1}{p^*}\intor b(x)|u_n|^{p^*}dx=c+ o(1)\ \mbox{as}\ n\ri\infty\end{array}\end{equation}
and
\begin{equation}\label{2}\begin{array}{ll}
\di\intor|\nabla u_n|^pdx&\di -\lambda \intor V(x)|u_n|^pdx-\intor a(x)|u_n|^rdx\\
&\di -\intor b(x)|u_n|^{p^*}dx= o(\|u_n\|)\ \mbox{as}\ n\ri\infty
\,.\end{array}\end{equation}
Relations \eqref{1} and \eqref{2} yield
\begin{equation}\label{3}\left(1-\frac pr\right)\intor a(x)|u_n|^rdx+\left(1-\frac{p}{r^*}\right)\intor b(x)|u_n|^{p^*}dx=O(1)+o(\|u_n\|)\ \mbox{as}\ n\ri\infty
\,.\end{equation}
Using hypothesis \eqref{hypr} in conjunction with the fact that the potentials $a$ and $b$ are positive, relation \eqref{3} implies
\begin{equation}\label{4}
\intor a(x)|u_n|^rdx=O(1)+o(\|u_n\|)\ \mbox{as}\ n\ri\infty\end{equation}
and
\begin{equation}\label{5}
\intor b(x)|u_n|^{p^*}dx=O(1)+o(\|u_n\|)\ \mbox{as}\ n\ri\infty
\,.\end{equation}

Inserting \eqref{4} and \eqref{5} in relation \eqref{2} we find
$$
\intor|\nabla u_n|^pdx-\lambda \intor V(x)|u_n|^pdx=O(1)+ o(\|u_n\|)\ \mbox{as}\ n\ri\infty
\,.$$
Now, since $\lambda<\lambda_1$ and using the minimization problem \eqref{mV}, we deduce that $(u_n)$ is bounded in $\d1pr$. Thus, up to a subsequence, we can assume that $(u_n)$ weakly converges  to some $u$ in $\d1pr$ and in $ L^{p^*}(\RR^N)$,
$$|\nabla u_n|^{p-2}\nabla u_n\rightharpoonup T\quad\mbox{in}\ (L^{p'}(\RR^N))^N$$
 and, by hypothesis \eqref{hypr},
$$u_n\ri u\quad\mbox{in} \ L^p_{loc}(\RR^N)\ \mbox{and}\ L^r_{loc}(\RR^N).$$
Moreover, $T$ and $u$ satisfy
\begin{equation}\label{relt}
-\mbox{div}\, T=\lambda V(x)|u|^{p-2}u+a(x)|u|^{r-2}u+b(x)|u|^{p^*-2}u\quad\mbox{in}\ (\d1pr)'.\end{equation}

By lower semicontinuity we find
$$\lambda\intor V(x)|u_n|^pdx+\intor a(x)|u_n|^rdx\ri \lambda\intor V(x)|u|^pdx+\intor a(x)|u|^rdx=:A\quad\mbox{as}\ n\ri\infty.$$
Relation \eqref{1} and our hypothesis $0<\lambda<\lambda_1$ imply that $A\geq 0$. We claim that $A>0$, provided that $c>0$ is small enough. Indeed, we first observe that relation \eqref{2} yields
\begin{equation}\label{7}\intor |\nabla u_n|^pdx=\intor b(x)|u_n|^{p^*}dx+A+o(\|u_n\|)\quad\mbox{as}\ n\ri\infty.\end{equation}
But relation \eqref{1} in combination with our assumption $\lambda\in (0,\lambda_1)$ imply that
$$\ell :=\lim_{n\ri\infty}\intor |\nabla u_n|^pdx>0.$$
Arguing by contradiction and assuming that $A=0$, relation \eqref{7} yields
$$\intor b(x)|u_n|^{p^*}dx\ri\ell\quad\mbox{as}\ n\ri\infty.$$
Returning to \eqref{1} we find that $c=\ell/N$. On the other hand, using the definition of the best Sobolev constant $S$, we have
$$\begin{array}{ll}
\di \ell &\di \geq S\lim_{n\ri\infty}\left(\intor |u_n|^{p^*}dx\right)^{p/p^*}=S\lim_{n\ri\infty}\left(\intor |u_n|^{p^*}dx\right)^{(N-p)/N}\\
&\di \geq S\,\|b\|^{(p-N)/N}_{L^\infty(\RR^N)}\lim_{n\ri\infty}\left(\intor b(x)|u_n|^{p^*}dx\right)^{(N-p)/N}=S\,\|b\|^{(p-N)/N}_{L^\infty(\RR^N)}\ell^{(N-p)/N},\end{array}$$
hence
$$\ell\geq S^{N/p}\,\|b\|_{L^\infty(\RR^N)}^{(p-N)/p}.$$
Since $\ell =cN$, in order to yield a contradiction with our assumption $A=0$, it suffices to choose
$c\in (0,c_0)$, where
\begin{equation}\label{czero}c_0:=\frac{S^{N/p}}{N}\,\|b\|_{L^\infty(\RR^N)}^{(p-N)/p}\,.\end{equation}
Fixing $c\in (0,c_0)$ we have $A>0$. Thus for some $R>0$,
$$\lim_{n\ri\infty}\sup_{z\in\RR^N}\int_{B_R(z)}\left(\lambda V(x)|u_n|^p+a(x)|u_n|^r\right)dx>0.$$

We have already seen that $u_n\rightharpoonup u$ in $\d1pr$ and $u_n\ri u$ almost everywhere. Passing again to a subsequence, we can assume that
$|\nabla u_n-\nabla u|^p\rightharpoonup \mu$ in ${\mathcal M}^+(\RR^N)$, $V^+|u_n- u|^p\rightharpoonup \nu$ in ${\mathcal M}^+(\RR^N)$,
$|\nabla u_n|^p\rightharpoonup \tilde\mu$ in ${\mathcal M}^+(\RR^N)$, and $| u_n|^{p^*}\rightharpoonup \tilde\nu$ in ${\mathcal M}^+(\RR^N)$. Set
$$\mu_\infty:=\lim_{R\ri\infty}\limsup_{n\ri\infty}\int_{\RR^N\cap (|z|>R)}|\nabla u_n|^pdx$$
and
$$\nu_\infty:=\lim_{R\ri\infty}\limsup_{n\ri\infty}\int_{\RR^N\cap (|z|>R)}V|u_n|^pdx.$$
Then by Lemma 2.1 in \cite{smets},

\smallskip
(i) $\mu_\infty\geq S_{\infty,V}\cdot\nu_\infty$.

\smallskip
(ii) $\nu=\sum_{i\in I}\nu_i\delta_{x_i}$ for some $x_i\in \Sigma_V$, $\nu_i>0$, $\mu\geq\sum_{i\in I}\nu_iS_V^{x_i}\delta_{x_i}$, $\tilde\mu\geq |\nabla u|^p+\sum_{i\in I}\nu_iS_V^{x_i}\delta_{x_i}$, and $\tilde\nu=|u|^{p^*}+\sum_{j\in J}\alpha_j\delta_{x_j}$ with $\alpha_j>0$ ($I$ and $J$ are at most countable).

\smallskip
(iii) $\limsup_{n\ri\infty}\intor V(x)|u_n|^pdx=\intor V(x)|u|^pdx+\|\nu\|+\nu_\infty$.

\smallskip
(iv) $\limsup_{n\ri\infty}\intor |\nabla u_n|^pdx=\intor |\nabla u|^pdx+\|\mu\|+\mu_\infty$ if $p=2$ and\\
$\limsup_{n\ri\infty}\intor |\nabla u_n|^pdx\geq \intor |\nabla u|^pdx+S_{*,V}\|\nu\|+\mu_\infty$ otherwise.

\smallskip
Returning to relations \eqref{1} and \eqref{2}, we obtain
$$\begin{array}{ll}
&\di \frac 1p\intor |\nabla u|^pdx+\frac 1p\sum_{i\in I}\nu_iS_V^{x_i}\leq c+\frac\lambda p\intor V(x)|u|^pdx+\frac\lambda p\,\|\nu\|+\frac\lambda p\,\nu_\infty+\\
&\di\frac 1r\intor a(x)|u|^rdx+\frac{1}{p^*}\intor b(x)|u|^{p^*}dx+\frac{1}{p^*}\sum_{j\in J}\alpha_jb(x_j)\end{array}$$
and
 \begin{equation}\label{sprefinal}\begin{array}{ll}
&\di \intor |\nabla u|^pdx+\sum_{i\in I}\nu_iS_V^{x_i}\leq \lambda \intor V(x)|u|^pdx+\lambda \,\|\nu\|+\lambda \,\nu_\infty+\\
&\di\intor a(x)|u|^rdx+\intor b(x)|u|^{p^*}dx+\sum_{j\in J}\alpha_jb(x_j).\end{array}\end{equation}
Combining these relations, we obtain
\begin{equation}\label{crelation}
\begin{array}{ll}c&\di\geq\frac 1N\intor b(x)|u|^{p^*}dx+\frac 1N\sum_{j\in J}\alpha_jb(x_j)+\left(\frac 1p-\frac 1r\right)\intor a(x)|u|^rdx\\
&\di\geq\frac 1N\intor b(x)|u|^{p^*}dx+\frac 1N\sum_{j\in J}\alpha_jb(x_j).\end{array}\end{equation}

Since $\ee'(u_n)\ri 0$ in $(\d1pr)'$ we deduce that for all $\phi\in C^\infty_0(\RR^N)$
$$\intor uT\cdot\nabla \phi dx+\intor\phi d\tilde\mu=\intor\phi bd\tilde\nu+\lambda\intor V(x)|u|^pdx.$$
Using now \eqref{relt} we obtain
$$\intor (uT\cdot\nabla\phi+\phi T\cdot\nabla u)dx=\lambda\intor V(x)|u|^pdx+\intor a(x)|u|^r\phi dx+\intor b(x)|u|^{p^*}\phi dx.$$
Combining these relations we find
\begin{equation}\label{almost}\begin{array}{ll}
\di  \intor\phi d\tilde\mu &=\di \intor \phi T\cdot\nabla udx-\intor b(x)|u|^{p^*}\phi dx+\intor\phi bd\tilde\nu\\
&\di\leq \intor \phi T\cdot\nabla udx+\intor\phi bd\tilde\nu .\end{array}\end{equation}
Concentrating $\phi$ on each $x_j$, relation \eqref{almost} yields $\nu_j\leq\alpha_jb(x_j)$. But for all $j$, we have $S\alpha_j^{p/p^*}\leq\nu_j$.
We deduce that
$$\alpha_j\geq S^{N/p}\left(b(x_j)\right)^{-N/p}\quad\mbox{for all}\ j\in J.$$
Thus if $J\not=\emptyset$, then relation \eqref{crelation} implies
$$c\geq \frac 1N\sum_{j\in J}\alpha_jb(x_j)\geq
\frac{S^{N/p}}{N}\|b\|_{L^\infty(\RR^N)}^{(p-N)/p},$$
which contradicts \eqref{czero} and the choice of $c\in (0,c_0)$. This shows that $J$ is empty, hence $\intor |u_n|^{p^*}dx\ri\intor |u|^{p^*}dx$. Using Proposition 3.32 from Brezis \cite{brebook} (which is a consequence of the Milman-Pettis theorem), we deduce that $u_n\ri u$ strongly in $L^{p^*}(\RR^N)$. We show that this implies the strong convergence of $(u_n)$ in $\d1pr$. For this purpose we employ an argument used in Filippuci, Pucci and R\u adulescu \cite[p. 713]{fpr}.
Consider the following elementary inequality (see formula (2.2) in Simon \cite{simons}): for all $\xi$, $\zeta\in\RR^N$
\begin{equation}\label{diaz1}
|\xi-\zeta|^{p}\leq\begin{cases} c(|\xi|^{p-2}\xi- |\zeta|^{p-2}\zeta
)(\xi-\zeta)\qquad&\mbox{for }\phantom{1<\,}p\geq 2;\\
c\langle|\xi|^{p-2}\xi-|\eta|^{p-2}\eta,\xi-\eta \rangle^{p/2}
\left(|\xi|^p+|\eta|^p\right)^{(2-p)/2}\qquad&\mbox{for }1<p<2,\end{cases}\end{equation}
where $c$ is a positive constant.

Restricting to the case $p\geq 2$, inequality \eqref{diaz1} implies that for all positive integers $n$ and $m$,
\begin{equation}\label{sd1p}\|u_n-u_m\|\leq |\ee'(u_n)(u_n-u_m)|+|\ee'(u_m)(u_n-u_m)|+|(\ee_0'(u_n)-\ee_0'(u_m))(u_n-u_m)|,\end{equation}
where $\ee_0:=\ee(u)-p^{-1}\intor |\nabla u|^pdx$. Applying the strong convergence of $(u_n)$ in $L^{p^*}(\RR^N)$, relation \eqref{sd1p} implies that $(u_n)$ strongly converges  in $\d1pr$. This concludes the proof of the Palais-Smale condition, provided that $c\in (0,c_0)$.\qed

Summarizing, in this section we have proved the following result.

\begin{lemma}\label{lemmaps}
Under the assumptions in Theorem \ref{t1}, the functional $\ee$ satisfies the Palais-Smale condition (PS)$_c$ for all $c\in (0,c_0)$, where $c_0=
\frac{S^{N/p}}{N}\,\|b\|_{L^\infty(\RR^N)}^{(p-N)/p}$.
\end{lemma}

Assuming that $1<p\leq N^2$ and following the same arguments as in the proof of Theorem 3.5 in Guedda and V\'eron \cite{guedda}, we can show that $\ee$ does not satisfy the localized Palais-Smale condition (PS)$_c$ if $c=\frac{kS^{N/p}}{N}\,\|b\|_{L^\infty(\RR^N)}^{(p-N)/p}$, for all positive integers $k$.

\section{Proof of the main result}
It remains to check the two geometric hypotheses of the mountain pass theorem. We have $\ee (0)=0$ and we argue the existence of a ``mountain" near the origin. For this purpose we first establish  that there are positive numbers $d$ and $r$ such that $\ee (u)\geq d$ for all $u\in \d1pr$ with $\|u\|=r$. Fix $0<\lambda<\lambda_1$. Using Theorem 3.1 from Smets \cite{smets}, there exists $\delta>0$ such that
\begin{equation}\label{prima}\intor |\nabla u|^pdx-\lambda\intor V(x)|u|^pdx\geq \delta\intor |\nabla u|^pdx\quad\mbox{for all $u\in\d1pr$}.\end{equation}
Taking into account the continuous embeddings of $\d1pr$ into $L^r(\RR^N)$ and $L^{p^*}(\RR^N)$ we obtain for all $u\in\d1pr$
$$\ee (u)\geq\frac\delta p\,\|u\|^p-C\left(\|u\|^r_{L^r(\RR^N)}+\|u\|^{p^*}_{L^{p^*}(\RR^N)}\right).$$
Using assumption \eqref{hypr} we deduce that $\ee (u)\geq d$ for all $u\in \d1pr$ with $\|u\|=r$, for some positive numbers $d$ and $r$.

The difficult part is to prove the existence of a ``valley" over the mountain. This will be achieved by using hypothesis \eqref{hypb}, which describes the decay of the potential $b$ near its maximum point in relationship with the critical nonlinear term. Let $\phi\not=0$ be an arbitrary function in $\d1pr$. Then
$$\begin{array}{ll}
\di \ee (t\phi)&\di =\frac{t^p}{p}\left(\intor |\nabla\phi|^pdx-\lambda \intor V(x)|\phi|^pdx\right)\\
&\di -t^p\left(\frac{t^{r-p}}{r}\intor a(x)|\phi|^rdx+\frac{t^{p^*-p}}{p^*}\intor b(x)|\phi|^{p^*}dx\right)<0,\end{array}$$
for  large enough $t>0$.

In order to ensure the localized Palais-Smale condition (PS)$_c$, it remains to show that the upper bounds of $\ee$ are in $(0,c_0)$, where $c_0$ is defined in \eqref{czero}. More precisely, if $u_\varepsilon$ achieves the minimum $S$ in problem \eqref{6} (recall that $u_\varepsilon$ is defined in \eqref{ueps}) then we prove that there exists $\varepsilon>0$ small enough such that
\begin{equation}\label{final!}
\sup_{t>0}\ee (tu_\varepsilon)<c_0:=\frac{S^{N/p}}{N}\,\|b\|_{L^\infty(\RR^N)}^{(p-N)/p}\,.\end{equation}

Fix $\varepsilon>0$. By invariance, we remark that
\begin{equation}\label{invariance}\intor |\nabla u_\ep|^pdx=\intor |\nabla u_1|^pdx\quad\mbox{and}\quad \intor b(x)u_\ep (x)^{p^*}dx=\intor b(\ep x)u_1 (x)^{p^*}dx.\end{equation}

As we have just observed, $\sup_{t>0}\ee (tu_\varepsilon)>0$ and this is achieved at some $t(\varepsilon)>0$. We claim that the family $\{t(\ep)\}_{\ep>0}$ is bounded from below by a positive constant. Indeed, combining $\ee'(t(\ep)u_\ep)(u_\ep)=0$ with relations \eqref{prima} and \eqref{invariance}, we obtain
$$t(\ep)^{p^*-p}\intor b(x)u_\ep^{p^*}dx+t(\ep)^{r-p}\intor a(x)u_\ep^{r}dxdx\geq\delta \intor |\nabla u_1|^pdx>0.$$
Using \eqref{hypr}, we deduce our claim. A straightforward computation shows that $\{t(\ep)\}_{\ep>0}$ is bounded from above. More precisely, our assumption \eqref{hypb} implies that there is some $R>0$ such that for all $\ep>0$
$$t(\ep)\leq\left(\frac{\di\intor |\nabla u_1|^pdx}{\di 2^{-1}b(0) \int_{B_R(0)} u_\ep (x)^{p^*}dx} \right)^{(N-p)/p^2}.$$

We control the behaviour of $\ee (t(\ep) u_\ep)=\sup_{t>0}\ee (tu_\varepsilon)$ by observing that
$$\ee (t(\ep) u_\ep)=\Phi_1(\ep)+\Phi_2(\ep)+\Phi_3(\ep),$$
where
$$\Phi_1(\ep)=\frac{t(\ep)^p}{p}\intor |\nabla u_1|^pdx-\frac{t(\ep)^{p^*}}{p^*}\, b(0)\intor u_1^{p^*}dx;$$
$$\Phi_2(\ep)=\frac{t(\ep)^{p^*}}{p^*}\, b(0)\intor u_1^{p^*}dx- \frac{t(\ep)^{p^*}}{p^*}\intor b(\ep x)u_1^{p^*}dx;$$
$$\Phi_3(x)=-\frac{\lambda t(\ep)^p}{p}\intor V(x)u_\ep^pdx-\frac{t(\ep)^r}{r}\intor a(x)u_\ep^rdx.$$

In what follows we prove that the growth of $\ee (t(\ep) u_\ep)$ is given by $\Phi_1$, while $\Phi_2$ and $\Phi_3$ tend to zero as $\ep\ri 0$.

Note that the mapping $(0,\infty)\ni s\longmapsto C_1s^p-C_2s^{p^*}$ (where $C_1$, $C_2$ are positive constants) admits a maximum for
$$s=\left(\frac{C_1(N-p)}{C_2N}\right)^{(N-p)/p^2}.$$
Returning to $\Phi_1$ we deduce that
$$\begin{array}{ll}
\di \Phi_1(\ep)&\di\leq\frac{1}{N}b(0)^{(p-N)/p}\left(\intor |\nabla u_1|^pdx\right)^{N/p}\left(\intor u_1^{p^*}dx\right)^{(p-N)/p}\\
&\di =\frac{S^{N/p}}{N}\,\|b\|_{L^\infty(\RR^N)}^{(p-N)/p}=c_0.\end{array}$$

It remains to establish the asymptotic decay of $\Phi_2$ and $\Phi_3$ as $\ep\ri 0$. Using hypothesis \eqref{hypb} we obtain, for some $C>0$ independent of $\ep$,
$$
 \Phi_2(\ep)\di\leq C\ep^\eta\int_{\RR^N}|x|^\eta\left(1+|x|^{p/(p-1)}\right)^{(p-N)/p},$$
which shows that
$$\Phi_2(\ep)\leq C\ep^\eta\quad\mbox{if}\ N\not=\frac{pr}{r+1-p}$$
and
$$\Phi_2(\ep)\leq C\ep^\eta\log\frac 1\ep\quad\mbox{if}\ N=\frac{pr}{r+1-p}.$$
A similar computation based on assumption \eqref{hypa} shows that
$$\Phi_3(\ep)\leq C\ep^\eta\quad\mbox{if}\ N\not=\frac{pr}{r+1-p}$$
and
$$\Phi_3(\ep)\leq C\ep^\eta\log\frac 1\ep\quad\mbox{if}\ N=\frac{pr}{r+1-p}.$$

Combining these estimates we obtain \eqref{final!}. This concludes the proof.\qed

\subsection{Final remarks}

Due to the singular behaviour of the indefinite potential $V$, we cannot improve the global regularity of the weak solution $u$. In the special case when $V$ is bounded (or away from its singularities, in the general case), Theorem 2.2 of Pucci and Servadei \cite{servadei} implies that $u\in L^\infty_{loc}(\RR^N)$. By Moser iteration, with the same arguments as in the proof of Theorem 1.1 in Filippucci, Pucci and R\u adulescu \cite{fpr}, this implies that $u\in C^{1,\alpha}(\RR^N\cap B_R)$, for some $\alpha=\alpha(R)\in (0,1)$. In such a case, $u\in L^m(\RR^N)$ for all $p^*<m<\infty$ and $\lim_{|x|\ri \infty}u(x)=0$, with the same ideas as in the proof of Lemma 2 in Yu \cite{yu}, which is based on Theorem 1 of Serrin \cite{james}.

We point out that an existence result in relationship with our Theorem \ref{t1} is proved in Theorem 3.1 of Guedda and V\'eron \cite{guedda} in the case of {\it bounded} domains, with only one perturbation term, and with {\it constant positive} potentials. In their case, a positive solution vanishing on the boundary is found, provided that $1<p^2\leq N$.

The result stated in Theorem \ref{t1} can be extended with similar arguments in the following three directions:

 (i) If the nonlinearity $|u|^{r-2}u$ is replaced by a more general function $g(x,u)$ with upper and lower bounds of the type $g_1(x)u^{r_1}$ and $g_2(x)u^{r_2}$ satisfying appropriate technical conditions;

 (ii) In the proof of the Palais-Smale condition (PS)$_c$, the fact that any bounded sequence in $\d1pr$ contains a strongly convergent subsequence can be proved under the stronger assumption that the subcritical term
$|u|^{r-2}u$ is replaced by an {\it almost critical} nonlinearity $h(x,u)$, in the sense that $h(x,u)=o(|u|^{p^*-1})$ as $|u|\ri\infty$, uniformly for $x\in\RR^N$. Next, with similar arguments, the conclusion of Theorem \ref{t1} follows.

(iii) The existence result established in Theorem \ref{t1} remains valid if problem \eqref{problem} is replaced with the following quasilinear singular problem
\begin{equation}\label{prfinal}
\left\{\begin{array}{ll}
&\di -\mbox{div}\, (|x|^{-ap}|\nabla u|^{p-2}\nabla u)-\mu \,\frac{|u|^{p-2}u}{|x|^{p(a+1)}} =\frac{|u|^{q-2}u}{|x|^{bq}}+\lambda f(x,u)\quad\mbox{in}\ \Omega\\
&\di u=0\quad\mbox{on}\ \partial\Omega\,,\end{array}\right.
\end{equation}
where $0\in\Omega\subset\RR^N$, $N\geq 3$, is a bounded domain and $1<p<N$, $a<N/p$, $a\leq b<a+1$, $\lambda$ is a positive parameter, $0\leq\mu<\bar\mu:=[(N-p)/p-a]^p$, $q=p^*(a,b):=Np/(N-pd)$ is the critical Hardy-Sobolev exponent and $d=a+1-b$. Note that $p^*(0,0)=p^*=Np/(N-p)$. In this case, $\lambda_1$ is the principal eigenvalue of the differential operator $L_\mu u:=-\mbox{div}\, (|x|^{-ap}|\nabla u|^{p-2}\nabla u)-\mu \,|x|^{-p(a+1)}|u|^{p-2}u$ and the role of the concentration-compactness principle of Smets \cite{smets} is played by Lemma 2.1 in Liang and Zhang \cite{liang}.

\smallskip
An interesting {\it open problem} is to study if the main result in the present paper remains true if
the $(p-1)$-superlinear term $|u|^{r-2}u$ is replaced by a nonlinear term $f(u)$ such that
$$\lim_{u\ri +\infty}\frac{f(u)}{u^{p-1}}=+\infty.$$

\medskip
{\bf Acknowledgements.} This research was supported by the Slovenian Research Agency grants P1-0292-0101, J1-5345-0101 and J1-6721-0101. We thank the referees for comments.

\vskip 3.0cm

\hfill 

Faculty of Education, University of Ljubljana, and Institute of Mathematics, Physics and Mechanics,
Kardeljeva plo\v{s}\v{c}ad 16, SI-1000 Ljubljana, Slovenia. 

Email: {\tt matija.cencelj@guest.arnes.si}\\

Faculty of Education and Faculty of  Mathematics and Physics,
University of Ljubljana, Kardeljeva plo\v{s}\v{c}ad 16, SI-1000 Ljubljana, Slovenia. 

Email: {\tt dusan.repovs@guest.arnes.si}\\

Faculty of Computer Science and Informatics, and Faculty of  Mathematics and Physics,
University of Ljubljana, Jadranska 21, SI-1000 Ljubljana, Slovenia. 

Email: {\tt ziga.virk@fmf.uni-lj.si}  

\end{document}